# Arbitrary Truncated Levy Flight: Asymmetrical Truncation and High-Order Correlations

**Dmitry V. Vinogradov**[1]

*Nizhny Novgorod Econophysics Laboratory,*

*11-96 Usilova st., Nizhny Novgorod, 603093, Russia*

**Abstract.**

The generalized correlation approach, which has been successfully used in statistical radio physics to describe non-Gaussian random processes, is proposed to describe stochastic financial processes. The generalized correlation approach has been used to describe a non-Gaussian random walk with independent, identically distributed increments in the general case, and high-order correlations have been investigated. The cumulants of an asymmetrically truncated Levy distribution have been found. The behaviors of asymmetrically truncated Levy flight, as a particular case of a random walk, are considered. It is shown that, in the Levy regime, high-order correlations between values of asymmetrically truncated Levy flight exist. The source of high-order correlations is the non-Gaussianity of the increments: the increment skewness generates threefold correlation, and the increment kurtosis generates fourfold correlation.

*Keywords:* Financial Stochastic processes, Truncated Levy Flights, High-Order Correlations, Non-Gaussian Random Walk

## 1. Introduction

A truncated Levy flight belongs to the class of discrete random walks with independent, identically distributed (i.i.d.) increments. The random walk is widely used as a description of physical and non-physical stochastic processes [1].

For example, the physical processes of Brownian motion and diffusion are described by the Gaussian random walk, which is a stochastically self-similar process with a fractal dimension equal to 2 [2].

Stochastically self-similar processes with other fractal dimensions occur every so often in non-physical systems. These are so-called Levy flights [3]. They have infinite variance, and their increments are distributed by the $\alpha$ - stable Levy laws of index $0 < \alpha < 2$ (Levy distributions).

The stochastic processes generated by economical and financial systems are the subject of investigation of the new interdisciplinary approach called econophysics [4, 5] and exhibit a number of particular features. The moments of their increments are finite, but the processes themselves are non-Gaussian. Their large-scale fluctuations are close to Brownian motion, whereas the small-scale fluctuations exhibit some of the Levy flight characteristics.

A model called truncated Levy flight has been proposed for the description of such stochastic processes for the first time in reference [6].

It should be noted that a number of financial stochastic process models are used in econophysics. The simplest Bachelier model is a Gaussian random walk (see, for example, [4]). In the 1960s Mandelbrot [7] and Fama [8] proposed a Levy flight model. In general, models must comply with the empirical behavior of common financial stochastic processes, called "stylized facts" [9, 10]. A truncated Levy flight is an idealized mathematical model and satisfies only some of the stylized facts, namely, "the absence of linear correlations between asset returns", "heavy tails" and "aggregational Gaussianity".

---

[1] Corresponding author. tel: +7(920)0757577, e-mail: Dmitry.Vinogradov@list.ru

The probability distribution of truncated Levy flight increments is a slightly deformed Levy distribution. This deformation must change the variance of the resulting distribution from the infinite to the finite, and consequently, according to the generalized central limit theorem, the resulting distribution belongs to the Gaussian basin of attraction. The chosen deformation must suppress the "tails" of the Levy distribution and cannot deform the central part. In the pioneering article [6], the abrupt truncation of the Levy distribution tails was used for this purpose.

Instead of the abrupt truncation employed in [6], a smooth exponential regression towards zero was introduced in [11]. This made it possible to derive an analytic expression for stochastic process characteristic function. The scale-invariant non-i.i.d. process, called scale-invariant truncated Levy process was introduced in [12].

A number of articles where other types of truncation are suggested have appeared later. Actually, there are a great number of deformation shapes that are solved in the posed problem. Therefore, a whole class of distributions that can be called "arbitrary truncated Levy distributions" exists.

Notwithstanding the fact that truncated Levy flights have received wide acceptance for the description of financial stochastic processes in econophysics, research focusing on the influence of the deformation shape on the stochastic process characteristics has only recently been carried out [13]. The problem of an arbitrary truncated Levy flight description using the cumulant approach has been solved for an approximation of a one-dimensional probability distribution function (pdf).

However, as was shown in [14], although it is used frequently in econophysics, the one-dimensional pdf stochastic financial process description characterizes only static behaviors and does not take into account dynamical properties including inner statistical relationships. This approach is insufficient for correctly describing stochastic financial processes and their models. Indeed, the presence of "volatility clustering" and "autocorrelation in absolute returns" in stochastic financial processes [9, 10] indicates the existence of high-order statistical relationships between increments.

Therefore, to correctly describe stochastic financial processes and to improve upon their models, it is necessary to use an approach that considers inner statistical relationships in financial processes.

This article is a continuation of article [13] and has two goals. The first goal is to show the possibilities and advantages of the generalized correlation function approach with regard to the description of a typical stochastic financial process model. It should be noted that the generalized correlation function approach, which takes into account high-order statistical relationships, has been successfully used in statistical radio physics for the non-Gaussian stochastic process description. The second goal is to investigate the behavior of asymmetrically truncated Levy flights, including high-order statistical relationships, based on this approach.

The rest of the article is organized as follows. Section 2 describes the generalized correlation function approach. The interdependences between the generalized correlation function approach and probability distribution approaches, namely, the one-dimensional pdf approach and the two-dimensional (joint) pdf approach, are presented in sections 2.1 and 2.2 to provide a better description of the generalized correlation function approach.

A complete statistical description of a random walk with i.i.d. non-Gaussian increments, considering all inner statistical relationships, is given in section 3.

The cumulants of an asymmetrically truncated Levy distribution are found in section 4, using the method suggested in [13]. The behavior of an asymmetrically truncated Levy flight is considered in section 5.

## 2. Generalized correlation approach as an alternative description of stochastic processes

It is well known (see, for example [2, 14, 15]) that the multivariate or joint probability distribution function, $W_N(x_1,t_1;x_2,t_2;\ldots;x_N,t_N)$, can be used to describe a stochastic process $X(t)$.

The simplest one-dimensional or one-point pdf, $W_1(x,t)$, characterizes a stochastic process at separate fixed points in time and does not allow the process to be dynamic. The two-dimensional or two-point joint pdf, $W_2(x_1,t_1;x_2,t_2)$, gives a more complete stochastic process description and allows for probability relationships between two values of the process at arbitrary points in time, $t_1$ and $t_2$.

However, only an infinite set of multivariate probability functions can provide complete information about a stochastic process:

$$W_1(x,t),\quad W_2(x_1,t_1;x_2,t_2),\ldots W_N(x_1,t_1;x_2,t_2;\ldots;x_N,t_N),\ldots, \tag{1}$$

or the multidimensional characteristic functions

$$\theta_1(q,t),\theta_2(q_1,t_1;q_2,t_2),\ldots \theta_N(q_1,t_1;q_2,t_2;\ldots;q_N,t_N)\ldots, \tag{2}$$

which are their Fourier transforms

$$\theta_N(q_1,t_1;q_2,t_2;\ldots;q_N,t_N) = \tag{3}$$
$$= \int_{-\infty}^{\infty}\ldots\int_{-\infty}^{\infty} W(x_1,t_1;x_2,t_2;\ldots x_N,t_N)\exp\{i(x_1q_1+x_2q_2+\ldots x_nq_n)\}\,dx_1dx_2\ldots dx_N.$$

The $N$-dimensional element of sets (1) or (2) alone characterizes the inner statistical relationships between $N$ values of a stochastic process at arbitrary points in time. It should be noted that the $N$-point term includes all of the information contained in the lower-order terms.

Frequently, as an alternative method of describing non-Gaussian stochastic processes, generalized correlation or cumulant functions are more convenient. This approach was first suggested by R.L. Stratonovich in [15] and has been intensively developed by the Nizhny Novgorod radio physics scientific school.

According to [15], the generalizing correlation or cumulant functions, $\mathrm{K}_N(t_1,t_2,\ldots t_N)$, by definition, are the derivatives of the logarithm of the characteristic functions:

$$\mathrm{K}_N(t_1,t_2,\ldots,t_N) = (-i)^N\left[\frac{\partial^N \ln\theta_N(q_1,t_1,\ldots,q_N,t_N)}{\partial q_1\partial q_2\ldots\partial q_N}\right]_{q=0}. \tag{4}$$

A stochastic process itself can be unambiguously and exhaustively described by an infinite set of cumulant functions:

$$\mathrm{K}_1(t_1),\quad \mathrm{K}_2(t_1,t_2),\quad \mathrm{K}_3(t_1,t_2,t_3),\ldots,\mathrm{K}_N(t_1,t_2,\ldots t_N),\ldots \tag{5}$$

Any characteristic of a stochastic process, particularly the characteristic functions or probability distributions, can be obtained from the generalizing correlation functions.

In contrast to the multivariate probability function description (1), each new term in set (5) contains only new information about the stochastic process and does not replicate the information contained in previous terms.

To clarify the physical meaning of the cumulant functions, we express them in terms of moment functions of the stochastic process [15]. Thus, the first-order cumulant function is simply the mean of the stochastic process:

$$\mathrm{K}_1(t) = \langle X(t)\rangle \equiv m(t), \tag{6}$$

where the single brackets $\langle\ldots\rangle$ represent a statistical average.

The second-order cumulant function

$$\mathrm{K}_2(t_1,t_2) = \langle X(t_1)X(t_2)\rangle - \langle X(t_1)\rangle\cdot\langle X(t_2)\rangle \equiv \langle X(t_1)X(t_2)\rangle - \mathrm{K}_1(t_1)\cdot\mathrm{K}_1(t_2) \tag{7}$$

is an ordinary autocorrelation function [16] describing the linear (first-order) statistical relationship between two values $X(t_1)$ and $X(t_2)$ of the stochastic process at different points in time, $t_1$ and $t_2$. It is necessary to stress that, in the case of time-point equality (degeneracy), the second-order cumulant function is the variance, $K_2(t,t) \equiv D(t) = \sigma^2(t)$, where $\sigma(t)$ is the standard deviation.

High-order cumulant functions or generalized correlation functions describe more complicated statistical relationships (correlations) between values of a stochastic process. Each cumulant function describes a statistical relationship at a specific order. Thus, the third-order cumulant function describes the relationship between three values of the stochastic process (threefold correlation or second-order statistical relationship) [15]:

$$K_3(t_1,t_2,t_3) = \langle X(t_1)X(t_2)X(t_3) \rangle - K_1(t_1) \cdot K_2(t_2,t_3) - K_1(t_2) \cdot K_2(t_1,t_3) - K_1(t_3) \cdot K_2(t_1,t_2) - \\ - K_1(t_1) \cdot K_1(t_2) \cdot K_1(t_3). \tag{8}$$

As follows from expression (8), the threefold correlation is non-zero only if the following conditions are met. First, three values of the stochastic process must be statistically dependent, and second, that dependence must be more complicated than the dependence of the twofold correlations (autocorrelations).

The fourfold correlation function $K_4(t_1,t_2,t_3,t_4)$ describes the statistical relationship between four values of the stochastic process (third-order statistical relationship), which in turn does not involve low-order correlations: neither threefold correlations nor twofold correlations. It should be noted that the physical significance of high-order correlations decreases as their order increases.

Dimensionless correlation coefficients [15] give the quantitative characteristics of the corresponding high-order statistical relationships:

$$R_j(t_1,t_2,...t_j) = \frac{K_j(t_1,t_2,...t_j)}{\sigma(t_1)\sigma(t_2)...\sigma(t_j)}. \tag{9}$$

The main advantages of the generalizing correlation function approach to stochastic process description are as follows:

1) Cumulant functions have a clear physical meaning and can be defined, to a certain extent, independently of each other, as the unique "normal coordinates" of the statistical description.

2) Cumulant functions permit a simple description of non-Gaussian stochastic processes.

3) A certain "good" probability distribution approximation represented by a real function, corresponding to a finite truncated set of cumulants, exists, whereas a non-singular function, corresponding to a finite set of moments, does not exist.

These advantages allow one to use cumulant functions as independent elements of a stochastic process statistical description.

## 2.1 Generalized correlation approach with a one-dimensional probability distribution function

Let us consider a one-dimensional pdf and a two-dimensional joint pdf of a stochastic process from the cumulant viewpoint to present the interdependence between the generalized correlation approach and the traditional multivariate pdf description.

As mentioned above, the one-dimensional function $W_1(x,t)$ or $\theta_1(q,t)$ gives the simplest description of a stochastic process and characterizes a stochastic process at separate fixed points in time.

As follows from definition (4), the set of completely degenerate (all points in time coincide, $t_1 = t_2 = ... = t_j = t$) cumulant functions corresponds to the one-point description of a stochastic process:

$$K_1(t), \ K_2(t^{[2]}), \ K_3(t^{[3]}),...,K_N(t^{[N]}),..., \tag{10}$$

which are the cumulants $K_j(t^{[j]}) \equiv k_j(t)$ of the one-dimensional pdf $W_1(x,t)$ or the coefficients in a Taylor expansion of the logarithm of the one-dimensional characteristic function [14, 17]:

$$\theta_1(q,t) = \exp\left[\sum_{j=1}^{\infty} \frac{k_j(t)}{j!}(iq)^j\right] \qquad (11)$$

As mentioned above, the first cumulant $k_1(t) \equiv m(t)$ is the mean of the stochastic process, and the second cumulant $k_2(t) \equiv D(t)$ is its variance.

The high-order cumulants, $k_j(t)$ (degenerate generalized correlation functions), or, more precisely, the dimensionless cumulant coefficients (standardized cumulants), $\Lambda_j(t) = k_j(t)/\sigma^j(t)$ where $j \geq 3$, define the shape of the one-dimensional pdf $W_1(x,t)$. We note that a normal distribution is characterized by two first-order cumulants (the mean and the variance), while the third- and fourth-order cumulant coefficients $\Lambda_3(t)$, $\Lambda_4(t)$ are the skewness coefficient and the kurtosis coefficient, respectively. The higher-order cumulant coefficients describe more complicated differences between this distribution and the normal distribution.

## 2.2 Generalized correlation approach with a two-dimensional joint probability distribution function

The description obtained from a two-dimensional joint pdf, $W_2(x_1,t_1;x_2,t_2)$, for a random process is more complicated and characterizes all statistical relationships between two values of a random process at arbitrary points in time, $t_1$ and $t_2$. The Taylor expansion of the two-dimensional characteristic function logarithm

$$\theta_2(q_1,t_1,q_2,t_2) = \exp\left[\sum_{s=1}^{\infty} \frac{i^s}{s!} \sum_{j=0}^{s} C_s^j K_s\left(t_1^{[s-j]}, t_2^{[j]}\right) q_1^{s-j} q_2^j\right] \qquad (12)$$

gives the following infinite set of cumulant functions corresponding to the two-point description:

$$\begin{array}{cccccc}
K_1(t_1) & K_1(t_2) & & & & \\
K_2(t_1^{[2]}) & K_2(t_1,t_2) & K_2(t_2^{[2]}) & & & \\
K_3(t_1^{[3]}) & K_3(t_1^{[2]},t_2) & K_3(t_1,t_2^{[2]}) & K_3(t_2^{[3]}) & & \\
K_4(t_1^{[4]}) & K_4(t_1^{[3]},t_2) & K_4(t_1^{[2]},t_2^{[2]}) & K_4(t_1,t_2^{[3]}) & K_4(t_2^{[4]}) & \\
K_5(t_1^{[5]}) & K_5(t_1^{[4]},t_2) & K_5(t_1^{[3]},t_2^{[2]}) & K_5(t_1^{[2]},t_2^{[3]}) & K_5(t_1,t_2^{[4]}) & K_5(t_2^{[5]})
\end{array} \qquad (13)$$

..........................................................................................

In set (13), the cumulant functions of equal order are placed together in the corresponding lines for clarity. The first line shows a pair of first-order cumulant functions, the second line contains three second-order functions, etc.

Set (13) involves cumulant functions of different types. The cumulant functions of the first type placed at the edges of the lines, $K_j(t_1^{[j]}) \equiv k_j(t_1)$ and $K_j(t_2^{[j]}) \equiv k_j(t_2)$, are completely degenerate. These functions correspond to the one-dimensional pdf $W_1(x_1,t_1)$ and $W_1(x_2,t_2)$, respectively. This is a result of the fact that a two-dimensional pdf involves all of the information contained in the one-dimensional pdf.

The cumulant functions of the second type describe the diversity of statistical relationships between two chosen values of the stochastic process. Thus, the second-order cumulant function (autocorrelation), $K_2(t_1,t_2)$, in the center of the second line is not degenerate and describes the simplest linear statistical relationship. More complicated statistical relationships between two values of the stochastic process are described by high-order cumulant functions, which are placed at the other lines. The threefold correlation (third line) is represented by two partially degenerate third-order cumulant functions, $K_3(t_1,t_2,t_2) \equiv K_3(t_1,t_2^{[2]})$ and $K_3(t_1,t_1,t_2) \equiv K_3(t_1^{[2]},t_2)$. The fourfold correlations, $K_4(t_1,t_2^{[3]})$, $K_4(t_1^{[2]},t_2^{[2]})$, and $K_4(t_1^{[3]},t_2)$, are placed in the fourth line and so on.

### 3. Generalized correlation approach with an i.i.d. random walk

Truncated Levy flights belong to the class of random walks with independent, identically distributed increments $\{x_i\}$ and discrete time:

$$X(n) = \sum_{i=1}^{n} x_i, \tag{14}$$

where $n$ is the step number, and the increments $\{x_i\}$ are specified by the set of cumulants $\kappa_j$.

It should be noted that a random walk is a Markovian non-stationary process, which is generally non-Gaussian. If the mean of the increments is equal to zero, $\kappa_1 = 0$, the random walk is also a martingale.

An i.i.d. random walk is the simplest model of stochastic financial processes used in econophysics. Thus, a drift-free Gaussian random walk or Brownian motion is used in the well-known Bachelier model of stochastic financial process. The Bachelier model satisfies (by virtue of the independency of increments) only one stylized fact [9], namely, the "absence of autocorrelation".

A more complicated non-Gaussian random walk has increments with a non-Gaussian probability distribution. Their high-frequency fluctuations are non-Gaussian, but their low-frequency fluctuations show Brownian motion behavior due to the central limit theorem (see, for example, [18]). These models satisfy an additional stylized fact, namely, "aggregational Gaussianity".

An arbitrary symmetrically truncated Levy flight (see [13]), as a particular case of a non-Gaussian random walk, satisfies yet another stylized fact, namely, "heavy tails" [9][2], due to its unique increment probability distribution.

An asymmetrically truncated Levy flight is considered in this article as a particular case of a non-Gaussian random walk; due to the asymmetry of the increment pdf and under specific conditions, it can satisfy the stylized fact of "gain/loss asymmetry" [9], bringing the number of modeled stylized facts up to four. It should be noted that an asymmetric Levy flight model have been successfully used in econophysics [19, 21, 20].

It should be noted that the next step in adapting a random walk model to the requirements of the stylized facts is the allowance of non-stationarity in the financial stochastic process, known as "intraday seasonality" [10], [14]. The usage of a non-Gaussian random walk with non-stationary increments would allow for this adaptation.

Let us use the generalized correlation approach to construct a complete statistical description of a random walk with independent, identically distributed increments that considers all internal statistical relationships.

It should be stressed that all results and conclusions obtained in section 2 can naturally be extended to the discrete time case by replacing the continuous time $t$ with a discrete time $n$. Thus, a discrete time random walk is unambiguously and exhaustively described by the set of cumulant functions $K_j(n_1,n_2,\ldots,n_j)$; see expression (5).

Let us use the cumulant bracket behavior [17] for the random walk cumulant functions, thus obtaining

---
[2] It should be noted that there is no consensus of opinion among the econophysical community as to the presence of "heavy tails" in financial stochastic processes. Thus, the Boston school accepts their presence (see, for example, [4]), while the Houston school does not, see [14].

$$\mathrm{K}_j(n_1, n_2, \ldots, n_j) \equiv \langle\langle X(n_1), X(n_2), \ldots X(n_j)\rangle\rangle, \tag{15}$$

where $\langle\langle \ldots \rangle\rangle$ is a cumulant bracket, namely, the linearity of the argument

$$\left\langle\left\langle \xi, \eta, \ldots, \sum_j a_j \zeta_j, \ldots, \omega \right\rangle\right\rangle = \sum_j a_j \langle\langle \xi, \eta, \ldots, \zeta_j, \ldots, \omega \rangle\rangle, \tag{16}$$

where $\xi, \eta, \zeta, \omega$ are random values, and $a$ is a deterministic value. We set the equality to zero:

$$\langle\langle \xi, \eta, \ldots, y, \ldots, \omega \rangle\rangle = 0 \tag{17}$$

in the case that one of the random values, for example, $y$, is statistically independent of the other values in its argument.

Let us represent each random walk value as part of the bracket argument in (15), except for $X(m)$, where $m = \min\{n_1, n_2, \ldots n_j\}$, using a sum of two statistically independent random values:

$$X(n_j) = X(m) + \sum_{l=m+1}^{n_j} x_l. \tag{18}$$

After substituting expression (18) into (15) and using the cumulant bracket behaviors (16) and (17), we obtain the desired expression for the cumulant functions:

$$\mathrm{K}_j(n_1, n_2, \ldots n_j) = \kappa_j \min\{n_1, n_2, \ldots n_j\}. \tag{19}$$

The cumulant functions (19) completely describe a random walk with independent, non-Gaussian, identically distributed increments (14).

It should be stressed that in the case of a Gaussian random walk or Brownian motion, a stochastic process is described by only the first two cumulant functions, $\mathrm{K}_1(n)$ and $\mathrm{K}_2(n_1, n_2)$, namely, the mean

$$\mathrm{K}_1(n) = m_1 \cdot n, \tag{20}$$

where the first increment cumulant $\kappa_1 = m_1$ is the drift coefficient and the autocorrelation function, and

$$\mathrm{K}_2(n_1, n_2) = \sigma_0^2 \min\{n_1, n_2\}, \tag{21}$$

where the second increment cumulant $\kappa_2 = \sigma_0^2$ is the variance, and $\sigma_0$ is the increment standard deviation.

## 3.1 One-dimensional probability distribution function of a random walk

As indicated above, the generalized correlation approach (5) to describing a stochastic process is complete, and any statistical characteristics of a stochastic process can be derived from it.

Thus, the one-dimensional probability function of a non-Gaussian random walk $W_1(x, n)$ (see section 2.1) is completely described by the set of cumulants $k_j(n) \equiv \mathrm{K}_j(n^{[j]})$, which possesses, according to expression (19), a linear dependence on the discrete time:

$$k_j(n) = \kappa_j n. \tag{22}$$

Consequently, the one-point pdf of a Gaussian random walk is characterized by the first two cumulants: a linear increasing mean (20) and the variance:

$$k_2(n) \equiv D(n) = \sigma_0^2 n, \tag{23}$$

where the diffusion coefficient $\sigma_0^2$ is the increment variance. It should be noted that (23) is the well-known diffusion law for Brownian motion (see, for example, [2]).

The shape of the one-dimensional pdf $W_1(x, n)$ of a non-Gaussian random walk is defined by the cumulants as well, or, more precisely, by the dimensionless cumulant coefficients (standardized cumulants) $\Lambda_j(n)$, where $j \geq 3$. It follows from (22) that

$$\Lambda_j(n) = \frac{\lambda_j}{n^{j/2-1}}, \tag{24}$$

where $\lambda_j = \kappa_j/\sigma_0^j$ is the *j*-th-order cumulant coefficient of the increment pdf. The faster the high-order cumulant coefficients $\Lambda_j(n)$ tend to zero, the higher their orders. When the number of steps tends to infinity, the high-order cumulants vanish, and the one-dimensional pdf $W_1(x,n)$ is described by the first two cumulants only or becomes Gaussian. Expression (24) is none other than the cumulant approach of the central limit theorem.

### 3.2 First-order statistical relationship of a random walk

The one-dimensional pdf approach for an i.i.d. random walk represented above characterizes a stochastic process at separate fixed points in time and does not allow the process to be dynamic. The inner statistical relationships of the random walk are described by generalized correlation functions. All of the relationship types can be separated by their order, which is characterized by the corresponding correlation functions.

Thus, the simplest linear statistical dependence between random walk values, termed the first-order relationship, corresponds to the autocorrelation (second-order cumulant) function $K_2(n_1, n_2)$. This statistical relationship exists in any type of random walk, both Gaussian and non-Gaussian. As follows from (9) and (19), the ordinary autocorrelation coefficient, which is the quantitative characteristic, equals

$$R_2(n_1, n_2) = \frac{\min\{n_1, n_2\}}{\sqrt{n_1 n_2}}. \tag{25}$$

The random walk autocorrelation coefficient, as evident from (25), is asymmetrical about the time interval by virtue of the non-stationarity of the process (the process itself, but not its increments). Indeed, by fixing one time point $n_1 = m$ in (25) and introducing the time interval $\tau = n_2 - n_1$, we find that the correlation coefficient dependence on the time interval for the "past" (when $\tau < 0$) $R_2(m,\tau) = \sqrt{1+\tau/m}$ differs from the dependence for the "future" $(\tau > 0)$ $R_2(m,\tau) = 1/\sqrt{1+\tau/m}$ (see Fig. 1).

The time interval in which the correlation between two random process values still persists is characterized by the correlation time $\tau_c$. Defining the correlation time as the curve half-width at the level of $R_2(m,\tau_c) = 0.5$, we find that the random walk correlation time for the "future" $\tau_c^+ = 3m$ is greater than the correlation time for the "past" $\tau_c^- = 3m/4$ by a factor of four.

We can see that the correlation time has a linear dependence on the chosen time point $m$. It should be noted that the random walk correlation time never vanishes. This means that any two random walk values are correlated to some extent. One exception is the statistical relationship between any value and the deterministic value at the initial point in time $X(0) = 0$.

### 3.3 High-order statistical relationships of a random walk

More complicated statistical relationships are inherent in non-Gaussian random walks, and they correspond to high-order generalized correlation functions. As follows from (19), the high-order correlation coefficients of a random walk, as quantitative characteristics of the high-order statistical relationships, equal

$$R_j(n_1, n_2, \ldots n_j) = \frac{\lambda_j \min\{n_1, n_2, \ldots n_j\}}{\sqrt{n_1 n_2 \ldots n_j}}. \tag{26}$$

We can see from (26) that the high-order statistical relationships are generated by the non-Gaussianity of the random walk increments. Moreover, each high-order increment cumulant corresponds to "own" statistical relationship. Thus, the second-order relationship or threefold correlation is generated by the skewness of the increment probability distribution, and the third-order relationship or fourfold correlation is due to the kurtosis of the increments, etc. The degree of high-order statistical coupling between the random walk values is determined by the appropriate increment cumulant coefficients.

By virtue of the non-stationarity of a random walk, all correlation coefficients are asymmetrical about the time interval. For example, the second-order statistical relationship or threefold correlation coefficient is

$$R_3(n_1, n_2, n_3) = \frac{\lambda_3 \min\{n_1, n_2, n_3\}}{\sqrt{n_1 n_2 n_3}}, \qquad (27)$$

where $\lambda_3$ is the kurtosis coefficient of the increment probability distribution function.

The threefold correlation coefficient peaks in the case of total degeneracy $(n_1 = n_2 = n_3 = n)$ and equals the one-dimensional pdf skewness coefficient of a random walk at the time point $n$: $R_3(n^{[3]}) = \Lambda_3(n) = \lambda_3/\sqrt{n}$. We note that any three random walk values, except the deterministic value at the origin, exhibit a threefold correlation in the presence of increment asymmetry.

Let us fix one point in time, $n_1 = m$, and introduce two time intervals, $\tau_2 = n_2 - n_1$ and $\tau_3 = n_3 - n_1$, between three stochastic process values and consider the time dependence of the threefold correlation coefficient. As with the ordinary autocorrelation, due to the non-stationarity of the process, a difference between the dependence for the "past" and the "future" exists.

There are three types of time dependence, dividing the domain $(\tau_2, \tau_3)$ into three regions (see Fig. 1). In region A, composed of the first quadrant, where $\tau_2 > 0, \tau_3 > 0$ (the "absolute future" quadrant), the threefold correlation coefficient time dependence is

$$R_3(m, \tau_2, \tau_3) = \frac{\Lambda_3(m)}{\sqrt{(1 + \tau_2/m)(1 + \tau_3/m)}}. \qquad (28)$$

By virtue of the time interval identity, the correlation coefficient is symmetrical about the axis $\tau_2 = \tau_3$. In region B, composed of the second quadrant, where $\tau_2 < 0, \tau_3 > 0$ (the "past" and "future" quadrant) and the upper half of the third quadrant, where $\tau_2 < 0, \tau_3 < 0$ (the "absolute past" quadrant), about the axis $\tau_2 = \tau_3$, the threefold correlation coefficient time dependence is

$$R_3(m, \tau_2, \tau_3) = \Lambda_3(m) \sqrt{\frac{1 + \tau_2/m}{1 + \tau_3/m}}. \qquad (29)$$

In region C, composed of the fourth quadrant, where $\tau_2 > 0, \tau_3 < 0$ (the "future and "past" quadrant) and the bottom half of the third quadrant, where $\tau_2 < 0, \tau_3 < 0$ (the "absolute past" quadrant), about the axis $\tau_2 = \tau_3$, the threefold correlation coefficient time dependence is

$$R_3(m, \tau_2, \tau_3) = \Lambda_3(n_1) \sqrt{\frac{1 + \tau_3/m}{1 + \tau_2/m}}. \qquad (30)$$

It is clear that the threefold correlation coefficient is asymmetrical about the past and the future.

The threefold correlation time of the random walk values is characterized by an isoline on the domain $(\tau_2, \tau_3)$, which corresponds to a correlation coefficient value equal to half of its maximal value, $R_3(m, \tau_2, \tau_3) = \Lambda_3(m)/2$. The isoline has three angles, with coordinates $(3m, 0)$, $(0, 3m)$ and $(-m, -m)$, and consists of two lines, specified by the equations $\tau_3 = 4\tau_2 + 3m$ and $\tau_3 = (\tau_2 - 3m)/4$, and one curve, specified by the equation $\tau_3 = 4m^2/(m + \tau_2) - m$.

## 4. Cumulants of arbitrary asymmetrically truncated Levy flight

The generalized correlation approach of an i.i.d. random walk has been carried out in the above sections. This approach allows us to describe the inner statistical relationships of a random walk at any order. The truncated Levy flight, as a simple econophysical model, belongs to the class of random walks. Let us use the generalized correlation approach for an asymmetrically truncated Levy flight description.

As mentioned above, it is necessary to have a set of increment cumulants for a complete statistical description of an i.i.d. random walk. Thus, let us find the cumulants of an arbitrary asymmetrically truncated Levy distribution.

Let us assume a non-symmetrical arbitrary truncated Levy distribution as a product of two functions:

$$P(x) = C P_L(x) g(x), \qquad (31)$$

where $P_L(x)$ is the symmetrical not drifted $\alpha$- stable Levy distribution [22,23], $g(x)$ is the non-symmetrical arbitrary deformation function, and $C$ is a normalizing constant. Let us assume that the deformation function equals a unit at the origin $g(0)=1$, and it tends to zero simultaneously with all of their derivatives faster than every power of $1/|x|$ when the argument tends to the infinity $|x| \to \infty$. Let us also assume that the characteristic spatial scale $l$ of the deformation function $g(x)$ is many times greater than the characteristic spatial scale $\gamma$ of the initial Levy distribution:

$$l \gg \gamma. \qquad (32)$$

Let us represent a non-symmetrical deformation function like a sum, $g(x) = g_{even}(x) + g_{odd}(x)$, of even, $g_{even}(x)$, and odd, $g_{odd}(x)$, functions, where

$$g_{even}(x) = \frac{g(x) + g(-x)}{2}, \qquad (33)$$

and

$$g_{odd}(x) = \frac{g(x) - g(-x)}{2}. \qquad (34)$$

Further without loss of generality the odd and even parts of deformation function will be considered separately.

To find of the set of distribution cumulants $\kappa_j$, let us use the cumulant expression as a function of moments $m_j$ [24]:

$$\kappa_j = j! \sum_{i=1}^{j} \sum \left(\frac{m_{p_1}}{p_1!}\right)^{\pi_1} \cdots \left(\frac{m_{p_i}}{p_i!}\right)^{\pi_i} \frac{(-1)^{\rho-1}(\rho-1)!}{\pi_1! \ldots \pi_i!}, \qquad (35)$$

where the second summation is taken over all positive $\pi$ and $\rho$, which obey the conditions:

$$p_1 \pi_1 + p_2 \pi_2 + \ldots + p_i \pi_i = j \qquad (36)$$

and

$$\pi_1 + \pi_2 + \ldots \pi_i = \rho. \qquad (37)$$

and the lower cumulants are

$$\begin{aligned}
\kappa_1 &= m_1 \\
\kappa_2 &= m_2 - m_1^2 \\
\kappa_3 &= m_3 - 3 m_1 m_2 + 2 m_1^3 \\
\kappa_4 &= m_4 - 3 m_2^2 - 4 m_1 m_3 + 12 m_1^2 m_2 - 6 m_1^4
\end{aligned} \qquad (38)$$

The moments of distribution (31), in turn, can be found from the value of characteristic function $\theta(q)$ derivatives at the origin [18]:

$$m_j = (-i)^j \left[ \frac{d^j \theta(q)}{dq^j} \right]_{q=0}. \tag{39}$$

where $\theta(q)$ is the Fourier transform of the pdf (31).

Because the truncated Levy distribution (31) is expressed as a product of two functions, based on the properties of Fourier transforms [25], the characteristic function $\theta(q)$ is the convolution of their Fourier transforms:

$$\theta(q) = \frac{C}{2\pi} \cdot \int_{-\infty}^{+\infty} \theta_L(q-\tau) \cdot G(\tau) \, d\tau. \tag{40}$$

where $G(q) = \mathcal{F}[g] = \mathcal{F}[g_{even}] + \mathcal{F}[g_{odd}]$ is a Fourier transform of the deformation function $g(x)$ where $\mathcal{F}[g_{even}]$ is a real function and $\mathcal{F}[g_{odd}]$ is an imaginary function and

$$\theta_L(q) = \exp(-\gamma^\alpha |q|^\alpha), \tag{41}$$

is a characteristic function of the symmetric non-shifted $\alpha$-stable Levy distribution (see, for example [18-23]), where $\alpha$ is the index of stability and $\gamma$ is the spatial scale.

Based on the properties of the characteristic function [18], namely, $\theta(0) = 1$ and $\theta(-q) = \theta^*(q)$, where the sign $*$ means the complex conjugate value, and the symmetry of the Levy distribution, one obtains from expressions (31), (39)-(41):

$$m_j = (-i)^j \frac{\int_{-\infty}^{+\infty} \theta_L^{(j)}(q) \cdot G(q) \, dq}{\int_{-\infty}^{+\infty} \theta_L(q) \cdot G(q) \, dq}. \tag{42}$$

Furthermore, we can use the presence of a small parameter $\varepsilon = (\gamma/l)^\alpha \ll 1$, namely, the ratio of the spatial scales of the Levy distribution and the deformation function, and find the moments of the truncated Levy distribution similar to an asymptotic expansion of the small parameter $\varepsilon$.

Let us use the Laplace method (see for example [26]). The Fourier transform $G(q)$ of the deformation function is concentrated in the small neighborhood of the origin in the form of sharp maxima and vanishes away from this region. The sharper this maxima is, the smaller the parameter $\varepsilon$. On the other hand, the value of the characteristic function of the Levy distribution $\theta_L(q)$ undergoes small changes in this neighborhood. As a result, the exact value of function $\theta_L(q)$ can be replaced by an asymptotic expansion in the neighborhood of the origin, and the solution must be found as a power of the small parameter $\varepsilon$.

Let us enter dimensionless coordinates $\varsigma = l \cdot q$, $\xi = x/l$ and expand the characteristic function $\theta_L(\varsigma)$ with the asymptotic series of powers of the small parameter $\varepsilon$ to find

$$\theta_L(\varsigma) = \exp(-\varepsilon |\varsigma|^\alpha) \approx 1 - \varepsilon |\varsigma|^\alpha + \varepsilon^2 \frac{|\varsigma|^{2\alpha}}{2} - \ldots \tag{43}$$

We change the variables in expression (42) to dimensionless ones and replaces the exact value of the characteristic function by an approximate one (43). Retaining terms of order up to $\varepsilon^1$ inclusive, one rewrites expression (42) for the moments of the truncated Levy distribution as

$$m_j = -(-i)^j l^j \varepsilon \cdot \alpha(\alpha-1)\ldots(\alpha-j+1) \cdot \begin{cases} \int_{-\infty}^{\infty} |\varsigma|^{\alpha-j} G(\varsigma) d\varsigma & \text{for even } j \\ \int_{-\infty}^{\infty} |\varsigma|^{\alpha-j} \operatorname{sgn}(\varsigma) G(\varsigma) d\varsigma & \text{for odd } j \end{cases}. \tag{44}$$

Taking into account that the inverse Fourier transform $\mathcal{F}^{-1}$ of a power function [27] is

$$\mathcal{F}^{-1}\left[|\varsigma|^{\beta}\right] = -\frac{1}{\pi} \sin \frac{\pi\beta}{2} \Gamma(\beta+1) |\xi|^{-\beta-1}, \tag{45}$$

and

$$\mathcal{F}^{-1}\left[|\varsigma|^{\beta} \operatorname{sgn}(\varsigma)\right] = -\frac{i}{\pi} \cos \frac{\pi\beta}{2} \Gamma(\beta+1) |\xi|^{-\beta-1} \operatorname{sgn}(\xi), \tag{46}$$

where $\Gamma(x)$ is the Gamma function, and converting from the frequency domain to the time domain expression (44), the result is

$$m_j = l^j \varepsilon \frac{2}{\pi} \Gamma(\alpha+1) \sin \frac{\pi\alpha}{2} \cdot \begin{cases} \int_0^{\infty} |\xi|^{j-1-\alpha} g_{even}(\xi) d\xi & \text{for even } j \\ \int_0^{\infty} |\xi|^{j-1-\alpha} g_{odd}(\xi) d\xi & \text{for odd } j \end{cases}. \tag{47}$$

Several expressions concerning Levy distribution and truncated Levy distribution in the frame of used approximation are presented in the Appendix.

It should be noted that orders of the truncated Levy distribution moment magnitudes are $m_j : l^j \varepsilon$. Therefore, retaining terms of order up to $\varepsilon^1$, one obtains from expressions (35)-(37)

$$\kappa_j = m_j = l^{j-\alpha} \gamma^{\alpha} A(\alpha) \cdot \mu_j(\alpha), \tag{48}$$

where

$$A(\alpha) = \frac{2}{\pi} \Gamma(\alpha+1) \sin \frac{\pi\alpha}{2} \tag{49}$$

and

$$\mu_j(\alpha) = \int_0^{\infty} \xi^{j-1-\alpha} g(\xi) d\xi, \quad \text{where} \quad \begin{matrix} g(\xi) \equiv g_{even}(\xi) & \text{for } j-even \\ g(\xi) \equiv g_{odd}(\xi) & \text{for } j-odd \end{matrix} \tag{50}$$

where the function $\mu_j(\alpha)$ describes the truncation shape influence on the cumulants. It should be noted that the function $\mu_j(\alpha)$ is the Mellin transform of the deformation function $g(\xi)$ [25].

As mentioned above, the high cumulant coefficients describe the differences between this distribution and the normal distribution. Based on expressions (48)-(50), one obtains the high cumulant coefficients ($j \geq 3$) of truncated Levy distribution:

$$\lambda_j = \left(\frac{l}{\gamma}\right)^{\alpha(j-2)/2} \frac{\mu_j(\alpha)}{A(\alpha) \mu_2^{j/2}(\alpha)}. \tag{51}$$

Given the results obtained, the distinguishing feature of an arbitrary truncated Levy distribution, from the cumulant approach point of view, is the specified orders of cumulant coefficients magnitudes: $\lambda_j : \varepsilon^{1-j/2}$, (for example $\lambda_4 : \varepsilon^{-1}$, $\lambda_6 : \varepsilon^{-2}$ etc.). This dependence is a "visiting card" of a truncated distribution, and any probability distribution belongs to the class of truncated Levy distributions if their cumulant coefficients fulfill these requirements.

It follows from the results derived that the cumulants of the truncated Levy flight (48) are directly dependent on the spatial scale $l$ of the deformation function as well as on the ratio of spatial scales $\varepsilon = (\gamma/l)^{\alpha}$, namely, $\kappa_j : l^j \varepsilon$. It should be noted that by virtue of the small size of this ratio $(\gamma/l)^{\alpha} \ll 1$ the cumulant dependence $\kappa_j(\alpha) : (\gamma/l)^{\alpha}$ upon the initial Levy distribution stability index $\alpha$ is strong, whereas expression (49) describes the weak dependence upon the stability index $\alpha$. The cumulant dependence upon the deformation function shape is described by expression (50).

## 4.1 Examples of an asymmetrically truncated Levy distribution

Let us consider two sample deformation functions to illustrate the results obtained above.

*1. Mantegna – Stanley truncation.* The truncated Levy flight was first proposed in article [6], and an abrupt truncation of the distribution "tails" was used. The deformation function corresponding to the asymmetrical truncation used,

$$g_{ms}(\xi) = \begin{cases} 1, & -\beta \leq \xi \leq 1 \\ 0, & \xi < -\beta \text{ or } \xi > 1 \end{cases}, \tag{52}$$

where $\beta$ is the asymmetrical coefficient, generates the influence functions

$$\mu_j(\alpha) = \frac{1}{j-\alpha} \begin{cases} \dfrac{1+\beta^{j-\alpha}}{2} & \text{for } j-\text{even} \\ \dfrac{1-\beta^{j-\alpha}}{2} & \text{for } j-\text{odd} \end{cases}. \tag{53}$$

All cumulants of the given truncated distribution can be obtained from expressions (48)-(50), and, for example, the variance is

$$\sigma^2 = l^2 \left(\frac{\gamma}{l}\right)^\alpha A(\alpha) \cdot \frac{1}{2-\alpha} \cdot \frac{1+\beta^{2-\alpha}}{2}. \tag{54}$$

In the case of symmetrical truncation, $\beta = 1$, Eq. (54) coincides with the results derived in articles [6] and [13]. The mean or drift coefficient generated by asymmetrical truncation is

$$\kappa_1 = l \left(\frac{\gamma}{l}\right)^\alpha A(\alpha) \cdot \frac{1}{1-\alpha} \cdot \frac{1-\beta^{1-\alpha}}{2}. \tag{55}$$

Let us assume that the truncation asymmetry is small, $\beta = 1 + \delta$, $\delta = 1$. With an approximation that considers terms of order up to $o(\delta)$, the even-order influence functions equal the influence functions of a symmetrical truncation ($\beta = 1$):

$$\mu_j(\alpha) = \frac{1}{j-\alpha}, \tag{56}$$

but the odd-order influence functions do not depend on the order $j$ or the stability index of the Levy distribution $\alpha$:

$$\mu_j(\alpha) = -\frac{\delta}{2}. \tag{57}$$

Therefore, for small asymmetry, the mean and the variance of the Mantegna-Stanely truncated distribution are

$$m_1 = -l \left(\frac{\gamma}{l}\right)^\alpha A(\alpha) \cdot \frac{\delta}{2}, \tag{58}$$

and

$$\sigma^2 = l^2 \left(\frac{\gamma}{l}\right)^\alpha A(\alpha) \cdot \frac{1}{2-\alpha}. \tag{59}$$

The skewness and kurtosis coefficients, which are important characteristics of a truncated Levy distribution, are

$$\lambda_3 = -\left(\frac{l}{\gamma}\right)^{\alpha/2} \frac{(2-\alpha)^{3/2}}{\sqrt{A(\alpha)}} \frac{\delta}{2} \tag{60}$$

and

$$\lambda_4 = \left(\frac{l}{\gamma}\right)^{\alpha} \frac{1}{A(\alpha)} \frac{(2-\alpha)^2}{(4-\alpha)}. \tag{61}$$

*2. Exponential truncation.* Another significant example is exponential suppression, and the appropriate asymmetrical deformation function is

$$g_e(\xi) = \begin{cases} \exp(-\xi) & \xi \geq 0 \\ \exp\left(\dfrac{\xi}{\beta}\right) & \xi < 0 \end{cases}, \tag{62}$$

which generates the influence function

$$\mu_j(\alpha) = \Gamma(j-\alpha) \begin{cases} \dfrac{1+\beta^{j-\alpha}}{2} & \text{for } j - \text{even} \\ \dfrac{1-\beta^{j-\alpha}}{2} & \text{for } j - \text{odd} \end{cases}. \tag{63}$$

Correspondingly, in the case of exponential truncation, the variance and the mean are

$$\sigma^2 = l^2 \left(\frac{\gamma}{l}\right)^{\alpha} A(\alpha) \cdot \Gamma(2-\alpha) \frac{1+\beta^{2-\alpha}}{2} \tag{64}$$

and

$$\kappa_1 = l \left(\frac{\gamma}{l}\right)^{\alpha} A(\alpha) \cdot \Gamma(1-\alpha) \cdot \frac{1-\beta^{1-\alpha}}{2}. \tag{65}$$

In the case of small asymmetry, as in the above example, the even-order influence functions coincide with those of the symmetrical case:

$$\mu_j(\alpha) = \Gamma(j-\alpha). \tag{66}$$

In contrast, the odd-order influence functions depend on the order:

$$\mu_j(\alpha) = -\Gamma(j-\alpha+1) \frac{\delta}{2}. \tag{67}$$

The mean and the variance of the exponentially truncated distribution in the case of small asymmetry are

$$m_1 = -l \left(\frac{\gamma}{l}\right)^{\alpha} A(\alpha) \Gamma(2-\alpha) \cdot \frac{\delta}{2} \tag{68}$$

and

$$\sigma^2 = l^2 \left(\frac{\gamma}{l}\right)^{\alpha} A(\alpha) \cdot \Gamma(2-\alpha) \tag{69}$$

The skewness and kurtosis coefficients are

$$\lambda_3 = -\left(\frac{l}{\gamma}\right)^{\alpha/2} \frac{1}{\sqrt{A(\alpha)}} \frac{(2-\alpha)(3-\alpha)}{\Gamma(2-\alpha)^{1/2}} \frac{\delta}{2}, \tag{70}$$

and

$$\lambda_4 = \left(\frac{l}{\gamma}\right)^{\alpha} \frac{1}{A(\alpha)} \frac{(2-\alpha)(3-\alpha)}{\Gamma(2-\alpha)}. \tag{71}$$

## 5. Behavior of truncated Levy flight correlations

It is known [4-6] that behaviors of truncated Levy flight fluctuations depend on their scale $N$. When the large-scale fluctuations of a process have the nature of Brownian motion, it is called the *Gaussian regime* of truncated Levy flight [4,6]. The small-scale fluctuations have some properties of Levy flight fluctuations, and it is called the *Levy regime*.

***Gaussian regime.*** From the cumulant approach viewpoint, the Gaussian regime occurs when the high-order cumulant coefficients of the one-dimensional probability distribution $W_1(x,n)$ of the truncated Levy flight can be neglected: $\Lambda_j(n) = 1$.

As indicated above, the order of magnitude of the truncated Levy distribution cumulant coefficients is $\lambda_j : \varepsilon^{1-j/2}$, and the Gaussian regime, as described by (24), arises when the characteristic spatial scale $N_G$ of the fluctuations exceeds

$$N_G \gg (l/\gamma)^\alpha. \tag{72}$$

An asymmetrically truncated Levy flight in the Gaussian regime can be described by Brownian motion with drift using the two first cumulant functions: the mean, $K_1(n)$ (20), and the autocorrelation function, $K_2(n_1, n_2)$ (21). Correspondingly, the one-dimensional pdf in the Gaussian regime is specified by the two first cumulants: the linearly increased mean and the variance obeying the diffusion law (23). Let us recall that the drift coefficient is the first increment cumulant (the mean) and the diffusion coefficient is the second increment cumulant (the variance). Equations (54), (58-59), (64) and (68), (69) are the drift coefficients and variances in the particular cases of Mantegna-Stanely and exponential asymmetrical truncations.

It is interesting to note that the slightly asymmetrical ($\delta = 1$) truncated Cauchy flight (stability index $\alpha = 1$) has a mean equal to $m_1 = -\gamma\delta/2$ and a diffusion coefficient equal to $\sigma^2 = 2l\gamma/\pi$ both for Mantegna-Stanely truncation and exponential truncation.

Statistical relationships between truncated Levy flight values in the Gaussian regime are represented by only the simplest linear dependence (25) (autocorrelation). Due to the non-stationarity of the process, the first-order correlation time differs for the past and the future. For more details regarding the behaviors of autocorrelation, see section 3.1.

***Levy regime.*** In the case of small-scale fluctuations in truncated Levy flight, when the characteristic spatial scale of the fluctuations meets the condition $N \leq (l/\gamma)^\alpha$, the high-order cumulants of the pdf cannot be neglected.

In this case, it is necessary to employ the complete set of cumulant functions (19). We obtain the one-dimensional characteristic function of the truncated Levy flight from expressions (11), (22) and (48), taking into account condition (32):

$$\theta_1(q,n) = \exp\left[\sum_{j=1}^{\infty} \frac{n \cdot \gamma^\alpha l^{j-\alpha} A(\alpha) \cdot \mu_j(\alpha)}{j!}(iq)^j\right]. \tag{73}$$

Because the dependence $\kappa_j : l^{j-\alpha}\gamma^\alpha$ is typical for the truncated Levy distribution cumulants, we can obtain the following one-dimensional pdf from expression (73):

$$W_1(x,n) = P_L\left(\frac{x}{n^{1/\alpha}}\right) g(x), \tag{74}$$

where $P_L(x)$ is the undisturbed Levy distribution.

The physical meaning of equality (74) is as follows. In the Levy regime, the influence of the deformation function on the evolution of the Levy flight one-dimensional pdf is confined to modulations of its one-dimensional distribution. The time dependence of the truncated Levy one-dimensional pdf remains identical to that for ordinary Levy flight.

It should be noted that, under condition (72), the influence of the deformation function on the return is the second order of the smallness value, and the return equals that of an undisturbed Levy flight:

$$W_1(0,n) = \frac{\Gamma(1/\alpha)}{\pi\alpha\gamma n^{1/\alpha}}. \tag{75}$$

The statistical relationships between the values of truncated Levy flight in the Levy regime are more complicated than in the Gaussian regime. Apart from the autocorrelation, high-order correlations generated by the high-order increment cumulants exist in the Levy regime. For more details on the behavior of high-order correlations, see section 3.3.

As follows from (50), only asymmetrically truncated distributions have odd cumulants, and as a result, only asymmetrically truncated Levy flights possess odd-order correlations. Thus, the simplest odd high-order correlation is the threefold correlation (19), (27).

The magnitude of the threefold correlation depends on the increment skewness coefficient $\lambda_3$, which in turn depends on the shape of the deformation function. Thus, for the stochastic process of truncated Cauchy flight (stability index $\alpha = 1$), the maximal value of the threefold correlation coefficient for Mantegna-Stanely truncation, $R_3^{MS}\left(n^{[3]}\right) = -\sqrt{l/\gamma}\,\dfrac{\pi\delta}{4\sqrt{n}}$, is smaller than that for the exponential truncation by a factor of two, $R_3^{EX}\left(n^{[3]}\right) = -\sqrt{l/\gamma}\,\dfrac{\pi\delta}{2\sqrt{n}}$. Let us recall that the drift coefficients and the variances for these two cases are identical.

The fourfold correlation between values of truncated Levy flight (19), (27), as an even-order correlation, exists for any truncation (both symmetrical and asymmetrical). For the cases under consideration, the fourfold correlation for the Mantegna-Stanely truncation, $R_4^{MS}\left(n^{[4]}\right) = \pi l/6\gamma n$, is smaller than that for the exponential truncation by a factor of six (62), $R_4^{EX}\left(n^{[4]}\right) = \pi l/\gamma n$.

## 6. Conclusions

The generalized correlation approach, which has been successfully used in statistical radio physics for describing non-Gaussian random processes, has been proposed in this work for the description of stochastic financial processes. A brief description of this approach and its interconnection with other stochastic process descriptions, namely, approaches involving one-dimensional and two-dimensional probability functions, has been given.

Due to the advantages of the generalized correlation approach, it is possible to obtain a complete statistical description of a non-Gaussian random walk with independent, identically distributed increments. This approach describes all inner statistical relationships between random walk values.

It has been demonstrated that high-order statistical relationships between process values at different points in time are generated by the non-Gaussianity of the increment probability distribution. Moreover, each high-order relationship corresponds to a high-order increment cumulant. For example, the increment skewness generates a threefold correlation between the process values, and the increment kurtosis generates a fourfold correlation. The presence of correlation time asymmetry arising from the non-stationarity of a random walk has been shown.

A number of simple models of stochastic financial processes, including the arbitrarily truncated Levy flight, belong to the class of i.i.d. random walk processes. In article [13], the cumulants of symmetrical arbitrary truncated Levy flight were reported. Here, the cumulants of arbitrary asymmetrically truncated Levy flight have been obtained. A knowledge of these cumulants can allow for investigations of the behavior of asymmetrically truncated Levy flight, including high-order statistical relationships; this was carried out using the generalized correlation approach.

It has been shown that an asymmetrically truncated Levy flight in the Gaussian regime is a Gaussian random walk with drift; thus, the drift coefficient is defined by the mean of the increment distribution, and the diffusion coefficient is described by the increment variance. Particular cases of truncation, namely, Mantegna-Stanely truncation and exponential truncation, have been considered.

Using the cumulant approach, it has been shown that, in the Levy regime, the one-dimensional probability distribution function of truncated Levy flight has a sufficiently simple expression (74), with an approximation taking into account terms of order $\varepsilon^1$, where $\varepsilon = (\gamma/l)^\alpha \ll 1$, $\gamma$ is the characteristic spatial scale of the initial Levy distribution and $l$ is the characteristic spatial scale of the deformation

The physical meaning of equality (74) is as follows. In the Levy regime, the influence of the deformation function on the evolution of the Levy flight one-dimensional pdf is confined to modulations of its one-dimensional distribution. The time dependence of the truncated Levy one-dimensional pdf remains identical to that for ordinary Levy flight.

## 7. Appendix

Let us present several expressions concerning the Levy distribution and the truncated Levy distribution in the frame of an approximation (taking into account terms of order up to $\varepsilon^1$, where $\varepsilon = (\gamma/l)^\alpha \ll 1$ is a small parameter) and apply the technique of generalized functions.

The Levy distribution characteristic function is

$$\theta_L(\varsigma) \approx 1 - \varepsilon |\varsigma|^\alpha . \tag{A1}$$

The corresponding Levy probability distribution function is

$$l \cdot P_L(\xi) = \delta(\xi) + \frac{\varepsilon}{\pi} \sin\frac{\pi\alpha}{2} \Gamma(1+\alpha) |\xi|^{-\alpha-1}, \tag{A2}$$

where $\delta(\xi)$ is the delta function.

The characteristic function of the truncated Levy distribution is

$$\theta(\varsigma) \approx 1 + \frac{\varepsilon}{2\pi}\left(|\varsigma|^\alpha, G\right) - \frac{\varepsilon}{2\pi}|\varsigma|^\alpha * G, \tag{A3}$$

where the sign "$*$" denotes a convolution, and the corresponding probability distribution function is

$$l \cdot P(\xi) = \delta(\xi)(1-b) + \frac{\varepsilon}{\pi}\sin\frac{\pi\alpha}{2}\Gamma(1+\alpha)|\xi|^{-\alpha-1} \cdot g(\xi), \tag{A4}$$

where

$$b = \frac{\varepsilon}{\pi}\sin\frac{\pi\alpha}{2}\Gamma(1+\alpha)\left(|\xi|^{-\alpha-1}, g(\xi)\right). \tag{A5}$$

## 8. References

[1] Rudnick, J.A., Gaspari, G.D., Elements of the Random Walk: an Introduction for Advanced Students and Researchers, Cambridge University Press, Cambridge, 2004.
[2] Mazo, R.M., Brownian Motion: Fluctuations, Dynamics and Applications, Oxford University Press, Oxford, 2002.
[3] Chechkin, A.V., Gonchar, V.Y., Klafter, J., Metzler, R., Fundamentals of Levy Flight Processes, Advances in Chemical Physics, 133B (2006) 439.
[4] Mantegna, R.N., Stanley, H.E., An Introduction to Econophysics: Correlations and Complexity in Finance, Cambridge University Press, Cambridge, 2000.
[5] Bouchaud, J.-P., Potters M, Theory of Financial Risk and Derivative Pricing. Cambridge: University Press, Cambridge 2000.
[6] Mantegna, R.N., Stanley, H.E., Stochastic Process with Ultraslow Convergence to a Gaussian: The Trancated Levy Flight, Phys. Rev. Lett. 73 (1994) 2946.
[7] Mundelbrot, B., The variation of certain speculative prices, J. Bus. 36 (1963) 394
[8] Fama, E.F., The behavior of stock market prices, J. Bus. 38 (1965) 34
[9] Cont, R., Empirical properties of asset returns: stylized facts and statistical issues, Quantitative Finance 1 (2001) 2.
[10] Dacorogna, M.M., Gencay, R., Muller, U.A., Olsen, R.B., Pictet, O.V. An Introduction to High-Frequency Finance, Academic Press, San Diego, 2001.
[11] Koponen, I., Analytic Approach to the Problem of Convergence of Truncated Levy Flights Towards the Gaussian Stochastic Process, Phys. Rev. E 52 (1995) 1197.

[12] Podobnik, B., Ivanov, P.Ch, Lee, Y., Stanley, H.E., Scale-Invariant Truncated Levy Process, Europhys. Lett. 52 (2000) 491.
[13] Vinogradov, D.V., Cumulant Approach of Arbitrary Truncated Levy Flight, Physica A 389 (2010) 5794
[14] McCauley, J. L., Dynamics of Markets, Cambridge University Press, Cambridge, 2009 (second edition)
[15] Stratonovich, R.L., Topics in the theory of Random Noise, Gordon & Breach, NY, 1963
[16] Van Kampen N.G., Stochastic Processes in Physics and Chemistry, North-Holland, Amsterdam, 2007
[17] Gardiner, C.W., Handbook of Stochastic Methods for Physics, Chemistry and the Natural Sciences, Springer, Berlin, 2009.
[18] Feller, W., An Introduction to Probability Theory and Its Application, Vol. II, John Wiley and Sons, New York. 1971.
[19] Altman, E.I., Financial ratios, discriminant analysis and the prediction of corporate bankruptcy, J. Financ. 23 (1968) 589
[20] Zmijewski, M.E., Methodological issues related to the estimation of financial distress prediction models, J. Accounting Res. 22 (1984) 59
[21] Podobnik, B., Valentincic, A., Horvatic, D., Stanley, H.E., Asymmetric Levy flight in financial ratios, Proc. Natl. Acad. Sci. USA 108, (2011) 17883.
[22] Shiryaev, A.N., Essentials of Stochastic Finance: Facts, Models, Theory, World Scientific Publishing, Singapore, 1999.
[23] Zolotarev, V.M., One-Dimensional Stable Distributions, American Mathematical Society, 1992.
[24] Kendall, M.G., Stuart, A., The Advanced Theory of Statistics. V.1 Distribution Theory, Macmillan Pub Co., 1983.
[25] Korn, G.A., Korn, T.M., Mathematical Handbook, McGraw-Hill Book Company, 1968.
[26] Nayfeh, A.H., Introduction to Perturbation Techniques, John Wiley & Sons, 1981.
[27] Gelfand, I.M., Shilov, G.E., Generalized Functions, vol.1 Translated from the Russian, Academic Press, New York and London, 1964.

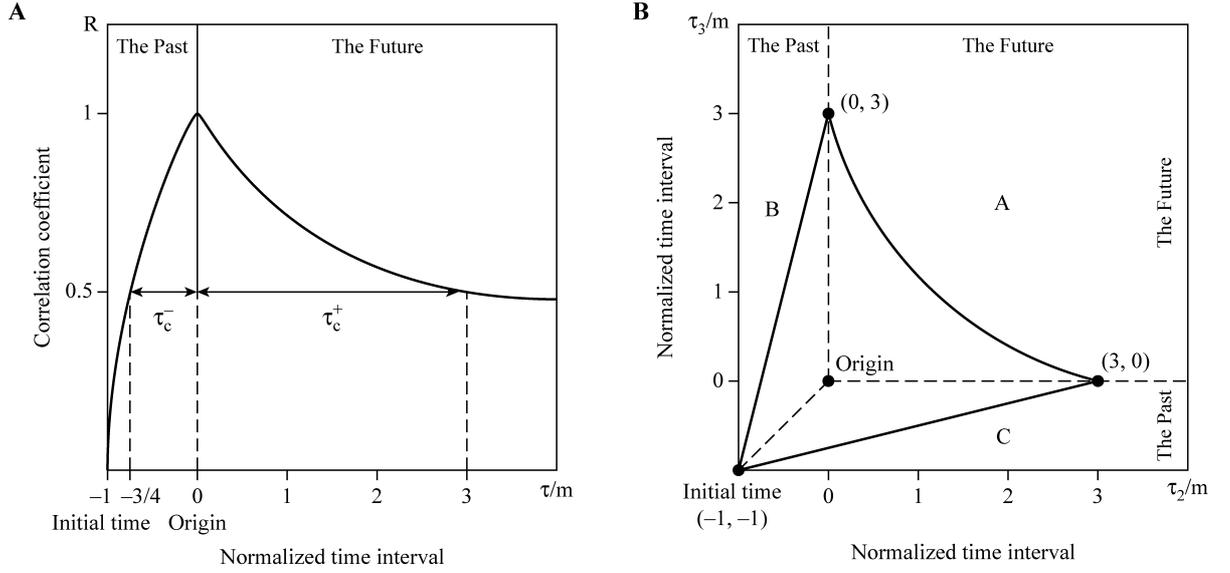

Fig. 1 Time dependence of non-Gaussian random walk correlations.

a) The time dependence of the first-order correlation (autocorrelation) coefficients of the random walk values for the time interval $\tau = n_2 - n_1$. The random walk correlation time for the "future" $\tau_c^+ = 3m$ is greater than the correlation time for the "past" $\tau_c^- = 3m/4$ by a factor of four.

b) Threefold correlations of non-Gaussian random walk values. There are three regions, A, B, and C, in the time interval $(\tau_2, \tau_3)$. The threefold correlation coefficient has different time dependences, see Eqs. (28), (29), and (30), in the different regions. The threefold correlation time is characterized by the isoline on the time interval domain, which corresponds to a correlation coefficient that is equal to half of its maximal value, $R_3(m, \tau_2, \tau_3) = \Lambda_3(m)/2$. The isoline has three angles with coordinates $(3m, 0)$, $(0, 3m)$, and $(-m, -m)$ and consists of two lines, specified by the equations $\tau_3 = 4\tau_2 + 3m$ and $\tau_3 = (\tau_2 - 3m)/4$, and one curve, specified by equation $\tau_3 = 4m^2/(m + \tau_2) - m$.